\documentclass[12pt]{amsart}
\usepackage{amssymb,latexsym,eufrak,amsmath,amscd, graphicx}
\usepackage[colorlinks=true,pagebackref,hyperindex]{hyperref}
\usepackage{amsfonts}
\usepackage{amscd}
\usepackage{xypic}

\renewcommand{\phi}{\varphi}
\renewcommand{\epsilon}{\varepsilon}
\renewcommand{\theta}{\vartheta}

\def\ZZ{{\mathbf Z}}

\def\AAA{{\mathbf A}}
\def\RR{{\mathbf R}}

\def\cJ{\mathcal{J}}

\def\cO{\mathcal{O}}

\def\fra{\mathfrak{a}}

\def\frq{\mathfrak{q}}

\DeclareMathOperator{\codim}{codim}

 \DeclareMathOperator{\lct}{lct}

\DeclareMathOperator{\Arn}{Arn}

\DeclareMathOperator{\ord}{ord}

\newtheorem{lemma}{Lemma}[section]
\newtheorem{theorem}[lemma]{Theorem}
\newtheorem{corollary}[lemma]{Corollary}
\newtheorem{proposition}[lemma]{Proposition}

\theoremstyle{definition}

\newtheorem{remark}[lemma]{Remark}

\newtheorem{example}[lemma]{Example}

\theoremstyle{remark}
\newtheorem*{remark*}{Remark}
\newtheorem*{note*}{Note}

\begin{document}

\title{A finiteness property of graded sequences of ideals}

\thanks{2000\,\emph{Mathematics Subject Classification}.
 Primary 14F18; Secondary 14B05. 
\newline 
The first author was partially supported by
NSF grants DMS-0449465 and DMS-1001740.
The second author was partially supported by
 NSF grant DMS-0758454 and
  a Packard Fellowship.}

\author[M.~Jonsson]{Mattias Jonsson}
\author[M.~Musta\c{t}\u{a}]{Mircea~Musta\c{t}\u{a}}
\address{Department of Mathematics, University of Michigan,
Ann Arbor, MI 48109, USA}
\email{mattiasj@umich.edu, mmustata@umich.edu}

\begin{abstract}
Given a graded sequence of ideals $(\fra_m)_{m\geq 1}$ on $X$, having
finite log canonical threshold, we show that if 
there are divisors $E_m$ over $X$ computing the log canonical threshold of $\fra_m$, and such that
the log discrepancies of the divisors $E_m$ are bounded, then the set $\{E_m\mid m\geq 1\}$ is finite.
\end{abstract}

\maketitle

\markboth{M.~JONSSON AND M.~MUSTA\c{T}\u{A}}{ON A FINITENESS PROPERTY
OF GRADED SEQUENCES OF IDEALS}

\section{Introduction}\label{intro}

Let $X$ be a smooth algebraic variety over an algebraically closed field $k$ of characteristic zero.
The log canonical threshold of a nonzero ideal $\fra$ on $X$ is a fundamental invariant of the singularities of the subscheme defined by $\fra$. Originally known as the \emph{complex 
singularity index}, 
it shows up in many contexts related to singularities, and it has found a plethora of applications in
birational geometry (see \cite{Kol} and \cite{EM}). 

In this note we will be interested in the behavior of this invariant in certain sequences of ideals.
Let $\fra_{\bullet}=(\fra_m)_{m\geq 1}$ be a \emph{graded sequence of ideals} on $X$, that is,
a sequence of ideals that satisfies $\fra_{\ell}\cdot\fra_m\subseteq\fra_{\ell+m}$ for every
$\ell, m\geq 1$. We always assume that, in addition, some ideal $\fra_m$ is nonzero. 
The main motivating example is the graded sequence $\fra^L_{\bullet}$ associated to a line bundle $L$ of nonnegative Iitaka dimension
on a smooth projective variety $X$: the ideal $\fra^L_m$ defines the base-locus of the
linear system $|L^m|$. Note that in this case the behavior of  $\fra^L_{\bullet}$
is easy to understand if the section ring $\oplus_{m}\Gamma(X,L^m)$ is finitely generated over $k$. 
Indeed, in this case there is a positive integer $p$ such that $\fra_{mp}=\fra_p^m$ for all $m$. The study of $\fra_{\bullet}^L$ is useful precisely when the section ring is not finitely generated
(or at least, when this property is not known \emph{a priori}). 

To a graded sequence $\fra_{\bullet}$ as above, one can associate an asymptotic version
of the log canonical threshold, by putting
$$\lct(\fra_{\bullet}):=\sup_{m;\fra_m\neq (0)}m\cdot\lct(\fra_m).$$
This can be infinite: for example, if $\fra_{\bullet}=\fra_{\bullet}^L$ as above, with $L$ big,
then $\lct(\fra_{\bullet})$ is infinite if and only if $L$ is nef (see Remark~\ref{nef} below). 

We will be concerned with the divisors that compute the log canonical thresholds
of the elements of a graded sequence. We denote by $A(\ord_E)$ the log discrepancy of 
a divisor $E$ over $X$ (see \S\ref{jumping_numbers} for the relevant definitions). The following is our main result, that gives
a positive answer to a question of Mihai P\u{a}un.\footnote{M.~P\u{a}un's question was motivated by the article \cite{Siu}, in which 
Y.-T.~Siu presents part of his arguments for the finite generation of the canonical ring.
At the end of \S 6.3, he evokes a subtle point in his approach, involving the control of an infinite
sequence of blow-ups. Although expressed in a different language, our main result shows that
the infinite blow-up process in Siu's approach can be ``stopped", provided that the log discrepancies 
of the divisors computing the log canonical thresholds are bounded.}

\noindent{\bf Theorem A}.
Let $\fra_{\bullet}$ be a graded sequence of ideals on a smooth variety $X$ such that 
$\lct(\fra_{\bullet})<\infty$. If $I\subseteq\ZZ_{>0}$ is a subset 
such that for all $m\in I$ we have a divisor $E_m$ over $X$ that computes $\lct(\fra_m)$
such that $\{A(\ord_{E_m})\mid m\in I\}$ is bounded, then the set
$\{E_m\mid m\in I\}$ is finite.

\noindent{\bf Corollary B}. Under the hypothesis in Theorem A, suppose that the set $I$ 
is infinite. Then there is a divisor $E$ over $X$ that computes $\lct(\fra_m)$ for infinitely $m$.
In particular, $E$ computes $\lct(\fra_{\bullet})$.

In fact, since our proof will require replacing $X$ by a suitable blow-up, we will need to prove a stronger version of the above theorem, in which we replace the log canonical threshold by the 
possibly
higher jumping numbers, in the sense of \cite{ELSV}  (see Theorem~\ref{main} below for the precise statement).

Here is a sketch of the proof. Let $Z_m$ be the image of $E_m$ on $X$, and let $W$ be the Zariski closure of $\bigcup_{m\in I}Z_m$. We may assume that $W$ is irreducible, and we first show 
that since $\lct(\fra_{\bullet})<\infty$, the asymptotic order of vanishing $\ord_W(\fra_{\bullet})$
is positive. In particular, $W$ is a proper subset of $X$. If $W$ has codimension at least two in $X$, then blowing-up $X$ along $W$ decreases the log discrepancies of the divisors $E_m$, and since these are bounded above, we reduce to the case 
when $W$ is a hypersurface. In this case, we use the following result, which we believe is of independent interest.

\noindent{\bf Theorem C}.
Let $H$ be a hypersurface in $X$, and $\fra$ a nonzero ideal. Suppose that $E$ is a divisor 
over $X$ that computes $\lct(\fra)$. If the image $Z$ of $E$ on $X$ is a proper subset of $H$,
and if $H$ is smooth at the generic point of $Z$, then the following inequality holds
$$\ord_Z(\fra)\geq\ord_H(\fra)\cdot \left(1+\frac{\ord_E(I_Z)}{A(\ord_E)}\right),$$
where $I_Z$ is the ideal defining $Z$.

Of course, as we have already mentioned, we need in fact a version of this result that applies 
also to higher jumping numbers (see Theorem~\ref{key_ingredient} below for this more general version of the theorem). Using Theorem~C, we show that if there were infinitely many $Z_m$
 that were properly contained in $W$, then the ideals in $\fra_{\bullet}$ would vanish along $W$ more than they should. Therefore all but finitely many of the $E_m$ are equal to $W$ (note that at this point we are on some blow-up of our original variety). 

In the following section we review some basic facts about log canonical thresholds and  higher jumping numbers. The proofs of the stronger versions of Theorems~C and A are given in 
\S\ref{an_inequality}, and respectively, \S\ref{the_main_result}.

\subsection{Acknowledgment}
We are grateful to Mihai P\u{a}un for the question that led to our main result.

\section{Jumping numbers and valuations}\label{jumping_numbers}

In this section we recall some definitions and results concerning the
invariants of singularities that we will use, and set the notation for the
rest of the paper.  
We work over a fixed algebraically closed field $k$ of characteristic zero.
Let $X$ be a smooth variety over $k$ (in particular, we assume that $X$ is connected
and separated). All ideal sheaves on $X$ are assumed to be coherent. 

By a divisor $E$ over $X$ we mean a prime divisor on a normal variety
$Y$ that has a proper birational morphism $\pi\colon Y\to X$. This induces a discrete valuation
of the function field $K(Y)=K(X)$, that we denote by $\ord_E$. As usual, we identify two such divisors
if they induce the same valuation. In particular, it follows from Hironaka's theorem on resolution of singularities that we may assume that both $Y$ and $E$ are nonsingular.
If we denote by $K_{Y/X}$ the relative canonical divisor, then the \emph{log discrepancy} 
of $\ord_E$ is given by $A(\ord_E):=1+\ord_E(K_{Y/X})$. Note that this depends on the variety $X$,
and whenever the variety is not clear from the context, we will write $A_X(\ord_E)$.
The \emph{center} of $E$ on $X$ is the image $c_X(E):=\pi(E)$ of $E$. 
We always consider on $c_X(E)$ the reduced scheme structure.
If $\fra$ is a nonzero ideal sheaf on $X$, we put
$$\ord_E(\fra):=\min\{\ord_E(f)\mid f\in\fra\cdot \cO_{X,c_X(E)}\}\in\RR_{\geq 0}.$$
If $Z$ is the subscheme defined by $\fra$, we also denote this by $\ord_E(Z)$.

Given an irreducible closed subset $Z$ of $X$, we define the order of vanishing along $Z$
as follows. Consider the normalized blow-up of $X$ along $Z$, and put $\ord_Z:=
\ord_E$, where $E$ is the unique irreducible component of the exceptional divisor that dominates 
$Z$. It is clear that in this case $c_X(E)=Z$. Note also that $\ord_Z(\fra)=\min_{x\in Z}\ord_x(\fra)$.

Let us recall the definition of multiplier ideals.
 For details and proofs we refer to
\cite[\S 9]{positivity}. Suppose that $\fra$ is a nonzero ideal on $X$. 
Let $\mu\colon X'\to X$ be a log resolution of $(X,\fra)$, that is, $\pi$ is proper and birational,
$X'$ is nonsingular, $\fra\cdot\cO_{X'}=\cO_{X'}(-F)$ for an effective divisor $F$, and 
$F+K_{X'/X}$ has simple normal crossings. For every $\lambda\in\RR_{\geq 0}$, the multiplier
ideal of $\fra$ of exponent $\lambda$ is given by 
$$\cJ(\fra^{\lambda}):=\pi_*\cO_{X'}(K_{X'/X}-\lfloor\lambda F\rfloor).$$
The definition is independent of the choice of log resolution.

It is clear from the above definition that if $\lambda<\lambda'$, then 
$\cJ(\fra^{\lambda'})\subseteq\cJ(\fra^{\lambda})$. Furthermore, for every $\lambda$
there is $\epsilon>0$ such that $\cJ(\fra^{\lambda})=\cJ(\fra^t)$ for every $t\in 
[\lambda,\lambda+\epsilon]$. One says that $\lambda>0$ is a \emph{jumping number} of
$\fra$ if $\cJ(\fra^{\lambda})\neq\cJ(\fra^{\lambda'})$ for every $\lambda'<\lambda$. 
It follows from the definition that if we write $F=\sum_ia_iE_i$, then for every jumping number
$\lambda$ there is $i$ such that $\lambda a_i$ is an integer. In particular, the jumping numbers form
a discrete set of rational numbers. 

For basic properties of the jumping numbers and applications, we refer to \cite{ELSV}. 
The most important jumping number is the smallest one, known as the
\emph{log canonical threshold} and denoted by $\lct(\fra)$. This is the smallest $\lambda$
such that $\cJ(\fra^{\lambda})\neq\cO_X$ (note that $\cJ(\fra^0)=\cO_X$). 

It is convenient to index the jumping numbers as follows (see \cite{JM}). 
Let $\frq$ be a nonzero ideal on $X$.
We put 
$$\lct^{\frq}(\fra):=\min\{\lambda\mid\frq\not\subseteq\cJ(\fra^{\lambda})\}.$$
Note that $\lct^{\cO_X}(\fra)$ is the log canonical threshold $\lct(\fra)$ of $\fra$. 
It follows from the definition that if $\fra\neq\cO_X$, then $\bigcap_{\lambda\geq 0}\cJ(\fra^{\lambda})
=(0)$, hence $\lct^{\frq}(\fra)$ is finite. When $\fra=\cO_X$, we make the convention
$\lct^{\frq}(\fra)=\infty$. We will also use the notation $\Arn^{\frq}(\fra):=1/\lct^{\frq}(\fra)$
(where $\Arn$ stands for \emph{Arnold multiplicity}).
It follows from the definition that we have
\begin{equation}\label{eq1_Section1}
\Arn^{\frq}(\fra)=\max_E\frac{\ord_E(\fra)}{A(\ord_E)+\ord_E(\frq)},
\end{equation}
where the maximum can be taken either over all divisors over $X$, or just over those lying
on a log resolution of $(X,\fra)$. We say that $E$ computes $\lct^{\frq}(\fra)$
(or $\Arn^{\frq}(\fra)$) if the maximum in (\ref{eq1_Section1}) is achieved by $E$.

The most interesting of the jumping numbers is the log canonical threshold. However, as the following lemma shows, the other jumping numbers appear naturally when we consider higher birational models.

\begin{proposition}\label{change_of_variable}
Let $\pi\colon X'\to X$ be a proper birational morphism, with $X'$ smooth, and $\fra$ and $\frq$ nonzero ideals on $X$. If $\fra'=\fra\cdot\cO_{X'}$, and
$\frq'=\frq\cdot\cO_{X'}(-K_{X'/X})$, then
$$\lct^{\frq}(\fra)=\lct^{\frq'}(\fra').$$
\end{proposition} 

\begin{proof}
This is an immediate consequence of (\ref{eq1_Section1}), and of the fact that for every divisor
$E$ over $X$, we have $A_X(\ord_E)=A_{X'}(\ord_E)+\ord_E(K_{X'/X})$.
\end{proof}

Suppose now that $\fra_{\bullet}$ is a graded sequence of ideals on $X$, and let 
$S=\{m\mid \fra_m\neq (0)\}$. Note that $S$ is closed under addition. 
In this case we have the following asymptotic version of the jumping numbers:
\begin{equation}\label{eq2_Section1}
\lct^{\frq}(\fra_{\bullet}):=\sup_{m\in S}m\cdot\lct^{\frq}(\fra_m)=
\lim_{m\to\infty, m\in S}m\cdot\lct^{\frq}(\fra_m)
\end{equation}
(see \cite[\S 2]{JM}). We put $\Arn^{\frq}(\fra_{\bullet})=1/\lct^{\frq}(\fra_{\bullet})$.
When $\frq=\cO_X$, we simply write $\lct(\fra_{\bullet})$ and $\Arn(\fra_{\bullet})$.
Note that $\lct^{\frq}(\fra_{\bullet})\in\RR_{>0}\cup\{\infty\}$. 
One can show that $\lct^{\frq}(\fra_{\bullet})=\infty$ if and only if $\lct(\fra_{\bullet})=\infty$
(see \cite[Corollary~6.10]{JM}). 

\begin{remark}\label{nef}
If $X$ is a smooth projective variety, $L$ is a big line bundle on $X$, and 
$\fra_{\bullet}=\fra_{\bullet}^L$ is the graded sequence of ideals defining the base loci of the powers of $L$ (see Introduction), then \cite[Corollary~2.10]{ELMNP} shows that
$\lct(\fra_{\bullet})=\infty$ if and only if
 $L$ is nef.
\end{remark}

If $\fra_{\bullet}$ is as above and $E$ is a divisor over $X$, we will also consider 
the following asymptotic version of the order of vanishing along $E$:
$$\ord_E(\fra_{\bullet}):=\inf_m\frac{\ord_E(\fra_m)}{m}=\lim_{m\to\infty,m\in S}\frac{\ord_E(\fra_m)}
{m}.$$
We have the following extension of (\ref{eq1_Section1})
\begin{equation}\label{Arn}
\Arn^{\frq}(\fra_{\bullet})=\sup_E\frac{\ord_E(\fra_{\bullet})}{A(\ord_E)+\ord_E(\frq)}.
\end{equation}
For these facts, we refer to \cite[\S 2]{JM}. We say that $E$ \emph{computes} 
$\lct^{\frq}(\fra_{\bullet})$ if the supremum in (\ref{Arn}) is achieved by $E$.
Note however that unlike in the case of one ideal, there may be no divisor $E$ that computes 
$\lct^{\frq}(\fra_{\bullet})$ (see \cite[Example~8.5]{JM}). 

We will use the following Izumi-type estimate (see \cite{izumi,ELS}).

\begin{proposition}\label{Izumi}
If $E$ is a divisor over $X$ with $c_X(E)=Z$, then
$$\ord_E(\fra)\leq A(\ord_E)\cdot \ord_Z(\fra)$$
for every nonzero ideal sheaf $\fra$ on $X$.
\end{proposition}

\begin{proof}
We may replace $X$ by an affine open subset of the generic point of $Z$, and therefore assume
that $X$ is affine. In this case we may assume that $\fra$ is principal.
If $\ord_Z(\fra)=m$, then for a general $p\in Z$ we have $\ord_p(\fra)=m$.
By \cite[Lemma~8.10]{Kol}, 
there is an open neighborhood $U$ of $p$ such that $\lct(\fra\vert_U)\geq 1/m$, and we get the assertion in the proposition since
$U\cap Z\neq\emptyset$ implies 
$\frac{A(\ord_E)}{\ord_E(\fra)}\geq\lct(\fra\vert_U)$. 
\end{proof}

\section{An inequality between orders of vanishing}\label{an_inequality}

We keep the notation and the conventions from \S\ref{jumping_numbers}.
The following is the main result in this section. Note that in the special case $\frq=\cO_X$,
this recovers Theorem~C in the Introduction. 

\begin{theorem}\label{key_ingredient}
Let $H$ be a hypersurface in $X$, and $\fra$, $\frq$ nonzero ideals on $X$. Suppose that $E$ is a divisor 
over $X$ that computes $\lct^{\frq}(\fra)$. If the center $Z$ of $E$ on $X$ is a proper subset of $H$,
and if $H$ is smooth at the generic point of $Z$, then the following inequality holds
\begin{equation}\label{formula_key_ingredient}
\ord_Z(\fra)\geq\ord_H(\fra)\cdot \left(1+\frac{\ord_E(Z)}{A(\ord_E)(1+\ord_H(\frq))}\right).
\end{equation}
\end{theorem}

We start by recalling a basic estimate for the log discrepancy of a valuation.
For a proof, see for example \cite[p. 157]{positivity}. 

\begin{lemma}\label{log_discrep}
Let $E$ be a divisor over $X$ with $c_X(E)=Z$, and let $\xi$ be the generic point of $Z$.
If $x_1,\ldots,x_r$ form a regular system of parameters of $\cO_{X,\xi}$, then
$$A(\ord_E)\geq\sum_{i=1}^r\ord_E(x_i).$$
\end{lemma}

\begin{corollary}\label{cor1}
If $H$ is a hypersurface in $X$, and  $E$ is a divisor over $X$ such that
$Z:=c_X(E)$ is a proper subset of $H$, and $H$ is smooth at the generic point of $Z$, then
\begin{equation}\label{eq_cor1}
A(\ord_E)\geq\ord_E(H)+\ord_E(Z).
\end{equation}
\end{corollary}

\begin{proof}
Let $\xi$ be the generic point of $Z$. Since $H$ is smooth at $\xi$, we may choose a regular 
system of parameters $x_1,\ldots,x_r$ of $\cO_{X,\xi}$ such that $H$ is defined at $\xi$ by $(x_1)$.
Note that by assumption $r\geq 2$.
By definition, we have $\ord_E(Z)=\min_j\ord_E(x_j)$. 
Let $i$ be such that $\ord_E(x_i)=\ord_E(Z)$. If $i\geq 2$, then by the lemma
$$A(\ord_E)\geq\ord_E(x_1)+\ord_E(x_i)=\ord_E(H)+\ord_E(Z).$$
On the other hand, if $i=1$, then using again the lemma we get
$$A(\ord_E)\geq \ord_E(x_1)+\ord_E(x_2)\geq 2\cdot\ord_E(x_1)=\ord_E(H)+\ord_E(Z).$$
\end{proof}

\begin{proof}[Proof of Theorem~\ref{key_ingredient}]
Let us put $m=\ord_H(\fra)$ and $p=\ord_H(\frq)$. We can write
$\fra=\cO_X(-H)^m\cdot \widetilde{\fra}$, and we get
\begin{equation}\label{eq0}
\ord_E(\fra)=m\cdot \ord_E(H)+\ord_E(\widetilde{\fra}),\,\,\,
\ord_Z(\fra)=m+\ord_Z(\widetilde{\fra})
\end{equation}
(note that $\ord_Z(H)=1$ since $H$ is smooth at the generic point of $Z$).
Since $E$ computes $\lct^{\frq}(\fra)$, it follows from (\ref{eq1_Section1}) that 
\begin{equation}\label{eq2_key_ingredient}
\frac{\ord_E(\fra)}{A(\ord_E)+\ord_E(\frq)}\geq\frac{\ord_H(\fra)}{A(\ord_H)+\ord_H(\frq)}=\frac{m}{1+p}.
\end{equation}
Corollary~\ref{cor1} gives $\ord_E(H)\leq A(\ord_E)-\ord_E(Z)$, and combining this with
(\ref{eq2_key_ingredient}) we deduce
\begin{equation}\label{eq4_key_ingredient}
m\leq (1+p)\cdot\frac{\ord_E(\fra)}{A(\ord_E)+\ord_E(\frq)}\leq
\frac{\ord_E(\widetilde{\fra})+m(A(\ord_E)-\ord_E(Z))+p\cdot \ord_E(\fra)}{A(\ord_E)+\ord_E(\frq)}
\end{equation}
$$
=m+\frac{\ord_E(\widetilde{\fra})+p\cdot \ord_E(\fra)-m(\ord_E(\frq)+\ord_E(Z))}{A(\ord_E)+
\ord_E(\frq)}.
$$
Therefore $\ord_E(\widetilde{\fra})\geq m(\ord_E(\frq)+\ord_E(Z))-p\cdot \ord_E(\fra)$. Using one more time the first equation in (\ref{eq0}), this implies
\begin{equation}\label{eq5_key_ingredient}
(1+p)\cdot \ord_E(\widetilde{\fra})\geq m(\ord_E(\frq)+\ord_E(Z))-pm\cdot \ord_E(H).
\end{equation}
On the other hand,
by Proposition~\ref{Izumi} we have $\ord_E(\widetilde{\fra})\leq A(\ord_E)\cdot \ord_Z(\widetilde{\fra})$, 
while clearly $\ord_E(\frq)\geq p\cdot \ord_E(H)$. Putting these together with (\ref{eq5_key_ingredient}) gives
$$
(1+p)A(\ord_E)\cdot\ord_Z(\widetilde{\fra})
\geq (1+p)\cdot \ord_E(\widetilde{\fra})$$
$$\geq 
m(\ord_E(\frq)+\ord_E(Z))-pm\cdot \ord_E(H)\geq m\cdot \ord_E(Z).
$$
Combining this with the second equality in (\ref{eq0}), we obtain
$$\ord_Z(\fra)=m+\ord_Z(\widetilde{\fra})\geq m\cdot
\left( 1+\frac{\ord_E(Z)}{A(\ord_E) (1+p)}\right),
$$
which completes the proof of the theorem.
\end{proof}

\begin{remark}
In Theorem~\ref{key_ingredient} one can replace $\ord_E$ by any real valuation 
of $K(X)$, having center on $X$ and 
computing $\lct^{\frq}(\fra)$. The proof goes through if one uses the definition of $A(v)$
from \cite[\S 5]{JM}. In this case, the assertion in Lemma~\ref{log_discrep}
follows from \cite[Corollary~5.4]{JM}. 
\end{remark}

\begin{example}\label{Example_key_ingredient}
The inequality in Theorem~\ref{key_ingredient} is optimal, at least in an asymptotic sense.
Indeed, let us consider the ideal $\fra=x^m(x,y^{m+1})$ in $k[x,y]$, where $m$ is a positive integer. Since this is a monomial ideal,
one can use Howald's theorem \cite{Howald} to compute its log canonical threshold. It is easy to check that $\lct(\fra)=\frac{m+2}{(m+1)^2}$, and this log canonical threshold is computed
by the (toric) divisor $E$ over $X=\AAA^2$ such that 
$$\ord_E(\sum_{i,j\geq 0}c_{i,j}x^iy^j)=\min\{(m+1)i+j\mid c_{i,j}\neq 0\}.$$
Note that $A(\ord_E)=m+2$, and the center of $E$ on $X$ is the origin.
If we take $\frq=\cO_X$ and $H=(x=0)$, then
$$\frac{\ord_Z(\fra)}{\ord_H(\fra)\left(1+\frac{\ord_E(Z)}{A(\ord_E)}\right)}=
\frac{m+1}{m\left(1+\frac{1}{m+2}\right)}=\frac{(m+1)(m+2)}{m(m+3)},$$
and this converges to $1$ when $m$ goes to infinity.
\end{example}

\begin{remark}
Note that the right-hand side of the inequality (\ref{formula_key_ingredient}) is
bounded above by $\ord_H(\fra)\cdot\left(1+\frac{1}{\lct(I_Z)\cdot(1+\ord_H(\frq))}\right)$,
where $I_Z$ is the ideal defining $Z$.
One could ask whether this expression is $\leq\ord_Z(\fra)$, improving in this way the assertion
in Theorem~\ref{key_ingredient}. However, this is not the case: let us consider the special
case $m=3$ in Example~\ref{Example_key_ingredient}, that is, $\fra=x^3(x,y^4)$. With $\frq=\cO_X$
and $H=(x=0)$, we have
$\ord_Z(\fra)=4$, while 
$$\ord_H(\fra)\cdot\left(1+\frac{1}{\lct(I_Z)}\right)=3\left(1+\frac{1}{2}\right)=\frac{9}{2}>4.$$
\end{remark}

\section{The main result}\label{the_main_result}

In this section we prove the generalized version of Theorem~A in the
 Introduction. We work in the same setting as in \S\ref{jumping_numbers}.

\begin{theorem}\label{main}
Let $\fra_{\bullet}$ be a graded sequence of ideals on $X$, and $\frq$ a nonzero ideal on $X$ such that $\lct^{\frq}(\fra_{\bullet})<\infty$. If $I\subseteq\ZZ_{>0}$ is a subset 
such that for all $m\in I$ we have a divisor $E_m$ over $X$ that computes $\lct^{\frq}(\fra_m)$
such that $\{A(\ord_{E_m})\mid m\in I\}$ is bounded, then the set
$\{E_m\mid m\in I\}$ is finite.
\end{theorem}

\begin{corollary}
Under the same hypothesis as in Theorem~\ref{main}, suppose that the set $I$ is infinite.
Then there is a divisor $E$ over $X$ that computes $\lct^{\frq}(\fra_m)$ for infinitely many $m$.
In particular, $E$ computes $\lct^{\frq}(\fra_{\bullet})$.
\end{corollary}

\begin{proof}[Proof of Theorem~\ref{main}]
Note that the hypothesis implies, in particular, that $\fra_m$ is nonzero for every $m\in I$.
We assume that $I$ is an infinite set, that $E_i\neq E_j$ for all $i\neq j$ in $I$
and aim to derive a contradiction. Let $Z_m=c_X(E_m)$. 
We argue by induction on $M:=\max\{A(\ord_{E_i})\mid i\in I\}$. This is finite by assumption.
Note that $M$ is a positive integer, and $M=1$
if and only if all the $E_i$'s are divisors on $X$. 
At several stages in the proof we will replace $I$ by an infinite subset. Note that this can only decrease the value of $M$.

We start with the following 
lemma.

\begin{lemma}\label{lem1_main}
With the above notation, suppose that there is an infinite subset $J\subseteq I$
such that $W:=\overline{\cup_{j\in J}Z_j}$ is irreducible, and $Z_j\neq W$ for all $j\in J$.
In this case 
$$\ord_W(\fra_{\bullet})\geq\Arn(\fra_{\bullet})\geq\Arn^{\frq}(\fra_{\bullet})>0.$$
\end{lemma}

\begin{proof}
We only need to prove the first inequality. Let $C=\Arn(\fra_{\bullet})$, so that 
$\Arn(\fra_m)\geq Cm$ for every $m$. If $j\in J$, then by Proposition~\ref{Izumi}
we have $\Arn(\fra_j)\leq\ord_{Z_j}(\fra_j)$. 

We need to show that $\ord_W(\fra_m)\geq Cm$
for every $m\geq 1$. We may, of course, assume that $\fra_m$ is nonzero. 
By hypothesis, we can find $0\leq\ell\leq m-1$ such that the set
\begin{equation}\label{set}
\bigcup_{j\in J, j\equiv \ell ({\rm mod}\,m)}Z_j
\end{equation}
is dense in $W$. Since all $Z_j$ are proper subsets of $W$, this implies that 
if in (\ref{set}) we only take the union over those $j\in J$ with $j\equiv\ell$ (mod $m$)
and with $j\geq N$, for some $N$, then the union is still dense in $W$.
 Let us fix $j_0\in J$ with $j_0\equiv \ell$ (mod $m$),
and let $C':=\max_{x\in W}\ord_x(\fra_{j_0})<\infty$ (recall that $\fra_{j_0}$ is nonzero).
If $mp+j_0\in J$, then the inclusion $\fra_m^p\cdot\fra_{j_0}\subseteq
\fra_{mp+j_0}$ implies
$$p\cdot \ord_{Z_{mp+j_0}}(\fra_m)+\ord_{Z_{mp+j_0}}(\fra_{j_0})\geq
\ord_{Z_{mp+j_0}}(\fra_{mp+j_0})\geq\Arn(\fra_{mp+j_0})\geq C(mp+j_0).$$
Therefore $\ord_{Z_{mp+j_0}}(\fra_m)\geq Cm-\frac{C'}{p}$.
Since we have arbitrarily large such $p$, and since the union of the corresponding $Z_{mp+j_0}$
is dense in $W$, we conclude that $\ord_W(\fra_m)\geq Cm$, as required.
\end{proof}

A first consequence of the lemma is that if $W$ is the closure of $\cup_{i\in I}Z_i$, then
$W\neq X$. In particular, this shows that when $M=1$, we have a contradiction. 

Arguing by Noetherian induction on $W$, we may assume that $W$ is minimal in $X$
with the property that there is an infinite family of divisors $(E_i)_{i\in I}$ as above, with 
$\max\{A(\ord_{E_i})\mid i\in I\}\leq M$. This implies first that
$W$ is irreducible. Indeed, if we 
consider the irreducible decomposition $W=W_1\cup\ldots\cup W_r$,
then there is $j$ such that $Z_i\subseteq W_j$ for infinitely many $i\in I$. Since we may
replace $I$
by $\{i\in I\mid Z_i\subseteq W_j\}$, it follows from the minimality assumption on $W$ that
$W=W_j$.

A second consequence of the minimality of $W$ is that for every infinite subset
$J\subseteq I$, the union $\cup_{j\in J}Z_j$ is dense in $W$. In particular, if
$U$ is an open subset of $X$ that meets $W$, then there are infinitely many $i\in I$
such that $U$ meets $Z_i$ (and the union of these $Z_i\cap U$ is dense in $W\cap U$).
Therefore in order to deduce a contradiction we may replace $X$ by $U$ and each $\fra_m$ by its restriction to $U$. We may thus assume that $W$ is nonsingular.

We claim that the induction hypothesis on $M$ implies that $W$ is a hypersurface in $X$.
Indeed, suppose that $c=\codim(W,X)\geq 2$, and let $\pi\colon X'\to X$ be the blow-up of $X$ along $W$. If $E$ is the exceptional divisor of $\pi$, then $K_{X'/X}=(c-1)E$.
Since $c_X(E_i)\subseteq W$ for every $i\in I$, it follows that $c_{X'}(E_i)\subseteq E$, hence
$$A_{X'}(\ord_{E_i})=A_X(\ord_{E_i})-\ord_{E_i}(K_{X'/X})\leq A_X(\ord_{E_i})-(c-1).$$
If $\fra'_m=\fra_m\cdot\cO_{X'}$ and $\frq'=\frq\cdot\cO_{X'}(-K_{X'/X})$, then by 
Proposition~\ref{change_of_variable} we have $\lct^{\frq}(\fra_i)=\lct^{\frq'}(\fra_i)$, and it follows from
hypothesis and (\ref{eq1_Section1}) that $E_i$ computes $\lct^{\frq'}(\fra'_i)$ for every $i\in I$.
Since $\max\{A_{X'}(\ord_{E_i})\mid i\in I\}\leq M-1$, we have a contradiction by induction on $M$. 

Therefore $W$ is a smooth hypersurface in $X$. If $Z_i=W$, then $E_i=W$, hence this can be the case for at most one $i$. After discarding this $i$, we may assume that each $Z_i$ is a proper subset of $W$. In particular, we may apply Theorem~\ref{key_ingredient} to get
\begin{equation}\label{eq11_main}
\ord_{Z_i}(\fra_i)\geq \ord_W(\fra_i)\cdot\left(1+\frac{\ord_{E_i}(Z_i)}{A(\ord_{E_i})(1+\ord_W(\frq))}\right).
\end{equation}

Note that $\ord_{E_i}(Z_i)\geq 1$ for all $i\in I$.
 Let $\alpha=\ord_W(\fra_{\bullet})$. We have $\alpha>0$ by Lemma~\ref{lem1_main}.
 Let us fix $\epsilon>0$ with $\epsilon<\frac{1}{M(1+\ord_W(\frq))}$. If we show that
 $\ord_W(\fra_m)\geq \alpha m(1+\epsilon)$ for every $m\geq 1$, then 
 $\alpha=\ord_W(\fra_{\bullet})\geq\alpha(1+\epsilon)$, a contradiction. 
 We now argue as in the proof of Lemma~\ref{lem1_main}. Let $0\leq\ell\leq m-1$ be such that the 
set in (\ref{set}) is dense in $W$. 
 We fix $j_0\in I$ such that $j_0\equiv\ell$ (mod $m$), and let $C':=\max_{x\in W}\ord_x(\fra_{j_0})$. It follows from the inclusion
 $\fra_{m}^p\cdot\fra_{j_0}\subseteq\fra_{mp+j_0}$ and from (\ref{eq11_main}) that 
 for every $p$ such that $mp+j_0\in I$ we have
 $$
 p\cdot \ord_{Z_{mp+j_0}}(\fra_m)\geq \ord_{Z_{mp+j_0}}(\fra_{mp+j_0})
 -\ord_{Z_{mp+j_0}}(\fra_{j_0})\geq \ord_W(\fra_{mp+j_0})(1+\epsilon)-C'.
 $$ 
 Therefore for every such $p$ we have $\ord_{Z_{mp+j_0}}(\fra_m)\geq \alpha m(1+\epsilon)-\frac{C'}{p}$.
 Since there are arbitrarily large such $p$, and the union
 of the corresponding $Z_{mp+j_0}$ is dense in $W$, we conclude that $\ord_W(\fra_m)\geq
 \alpha m(1+\epsilon)$. As we have seen, this leads to a contradiction, and thus completes the proof
 of the theorem.
\end{proof}

\providecommand{\bysame}{\leavevmode \hbox \o3em
{\hrulefill}\thinspace}


\begin{thebibliography}{positivity}


\bibitem[ELMNP]{ELMNP}
 L. Ein, R. Lazarsfeld, M.
Musta\c{t}\v{a}, M. Nakamaye and M. Popa, Asymptotic
invariants of base loci,  Ann.\ Inst. Fourier (Grenoble) \textbf{56} (2006), 
1701--1734.


\bibitem[ELS]{ELS}
L.~Ein, R.~Lazarsfeld, and K.~E.~Smith, Uniform approximation 
of Abhyankar valuations in smooth function fields, Amer. J. Math. \textbf{125} (2003),
409--440.



\bibitem[ELSV]{ELSV}
L.~Ein, R.~Lazarsfeld, K.~E.~Smith and D.~Varolin,
Jumping coefficients of multiplier ideals, Duke Math. J. \textbf{123}
  (2004), 469--506.


\bibitem[EM]{EM}
L.~Ein and M.~Musta\c{t}\u{a}, Invariants of singularities of pairs, in \emph{International Congress of Mathematicians},  Vol. II, 583--602, Eur. Math. Soc., Z\"{u}rich, 2006.


\bibitem[How]{Howald}
J.~Howald, Multiplier ideals of monomial ideals, Trans. Amer. Math. Soc. \textbf{353}
(2001), 2665--2671.

\bibitem[Izu]{izumi}
S.~Izumi,
A measure of integrity for local analytic algebras,
Publ. RIMS Kyoto Univ. \textbf{21} (1985), 719--735.


\bibitem[JM]{JM}
M.~Jonsson and M.~Musta\c{t}\u{a}, Valuations and asymptotic invariants for sequences of ideals,
arXiv: 1011.3699.

\bibitem[Kol]{Kol}
J.~Koll\'{a}r, Singularities of pairs, in \emph{Algebraic geometry,
Santa Cruz 1995}, 221--286,  Proc. Symp. Pure
Math. 62, Part 1, Amer. Math. Soc., Providence, RI, 1997.

\bibitem[Laz]{positivity}
R.~Lazarsfeld, \emph{Positivity in algebraic geometry} II, Ergebnisse der Mathematik und ihrer
Grenzgebiete, 3. Folge, Vol. \textbf{49}, Springer-Verlag, Berlin, 2004.

\bibitem[Siu]{Siu}
Y.-T.~Siu, Techniques for the analytic proof of the finite generation of the canonical ring, 
\emph{Current developments in mathematics, 2007}, 177--219, Int. Press, Somerville, MA, 2009.


\end{thebibliography}
\end{document}